\documentclass[12pt,a4paper]{article}
\usepackage{amssymb,amsmath,mathrsfs,amsthm,titling,changepage,graphicx,enumerate}
\usepackage[margin=0.5in]{geometry}
\usepackage[bottom]{footmisc}
\usepackage[affil-it]{authblk}

\usepackage{xcolor}
\newcommand*{\mybox}[2]{\colorbox{#1!30}{\parbox{.98\linewidth}{#2}}}

\newtheorem{lemma}{Lemma}

\newcommand{\subtitle}[1]{%
  \posttitle{%
    \par\end{center}
    \begin{center}\large#1\end{center}
    \vskip0.5em}%
}

\title{Zorn's Lemma}				
\subtitle{An elementary proof under the Axiom of Choice}
\author{Arjun Jain}	
\affil{\small{4th year Integrated MSc. Physics\\ Indian Institute of Technology Roorkee, Roorkee, India\\E-mail address: \texttt{arjunjain16@gmail.com}}}							
\date{}						

\begin{document}
\maketitle		

Here, I present an elementary proof of Zorn's Lemma under the Axiom of Choice, simplifying and supplying necessary details in the original proof by Paul R. Halmos in \emph{Naive Set Theory}\footnote{Paul R. Halmos. Naive Set theory. Van Nostrand Reinhold Company. 1960. Reprint by Martino Fine Books in  2011}.\newline\newline
To start with, I assume knowledge of basic Set Theory, i.e., axiom of extension, axiom of specification, axiom of pairing, axiom of unions, axiom of powers, ordered pairs, relations, functions, families, axiom of infinity, numbers, peano arithmetic, order, and the axiom of choice.\newline\newline
As the Axiom of Choice is central to the proof, here is the description, as given in Halmos's \emph{Naive Set Theory}:\newline\newline
\textbf{Axiom of Choice:}\newline
\emph{The Cartesian product of a non empty family of nonempty sets is non empty.}\newline\newline
Suppose that $\mathscr{C}$ is a non empty collection of non empty sets. We can convert $\mathscr{C}$ into an indexed 
set, by using the collection $\mathscr{C}$ itself in the role of the index set and using the identity mapping on $\mathscr{C}$ in the role of the indexing. The axiom of choice, then says that the Cartesian product of the sets of $\mathscr{C}$ has at least one element. An element of such a Cartesian product is, by definition, a function whose domain is the index set ($\mathscr{C}$) and whose value at each index belongs to the set bearing that index. 
Therefore, an equivalent form of the Axiom of Choice is that there exists a function f with domain $\mathscr{C}$ such that if $A\in\mathscr{C}$, then $f(A)\in A$. This conclusion applies, in particular, in case $\mathscr{C}$ is the collection of all non empty subsets of a non empty set X. The assertion in that case, is that there exists a function f with domain $P(X)\setminus\{\varnothing\}$ such that if A is in that domain, then $f(A)\in A$. In intuitive language, the function f can be described as a simultaneous choice of an element from each of many sets; this is the reason for the name of the axiom. 

\section{A Preamble to Zorn's Lemma}
\begin{adjustwidth}{2.5em}{0pt}
The Statement of Zorn's Lemma is as follows:\newline
\textbf{\emph{If X is a partially ordered set such that every chain in X has an upper bound, then X contains a maximal element.}}\newline\newline
Although Zorn's name has been stuck to this Lemma, there were, similar maximal principles(esp. by Hausdorff and Kuratowski), before Zorn published his results in 1933. But despite the fact, that there are many variants of Zorn's Lemma, they are all equivalent to each other and to the Axiom of Choice.\newline
So, how did the term ``Zorn's Lemma" come to be? Mycielski attributed to Semadeni\footnote{P.J. Campbell. The origin of "Zorn's lemma". Historia Math. , 5 (1978) pp. 77–89}, the following convincing explanation:\newline
\emph{Namely, in Science, the consumer decides upon the name of the tools which he uses, and the consumer is not always the best informed person.}\newline
Nonetheless, there is no doubt that Zorn provided a great service in directing attention to the largely unrealized potential in maximal principles.\newline\newline
As an example, Zorn's Lemma can be used to prove that every non zero vector space has a basis. 
We consider the set of all linearly independent subsets of the given vector space V, partially ordered by inclusion. Let Y be a chain of linearly independent subsets of V. We note that the union of such a set can serve as an upper bound for it. To apply Zorn's lemma, we have to check whether the union is linearly independent. If $t_1,\dotsc,t_n$ belong to the union, then each $t_i$ belongs to some linearly independent set $L_i\in Y$. Because Y is a chain, one of these sets $L_i$ contains all the others. If that is $L_j$, then the linear independence of $L_j$ implies that no non trivial linear combination of $t_1,\dotsc,t_n$ can be zero, which proves that the union of the sets in Y is linearly independent. Therefore, by Zorn's lemma, there is a maximal linearly independent set. Such a set is not just linearly independent, but also spans the whole space, since if it didn’t we could just pick an element that did not belong to its linear span and we could add it to the linearly independent set, contradicting maximality.
\end{adjustwidth}

\section{Proof of Zorn's Lemma:}

\begin{proof}\mbox{}\\*

First of all, the empty chain $\varnothing$, is a chain in X and by the hypothesis of Zorn's Lemma, must have an upper bound, say z, in X. As a result,  $X\neq\varnothing$ and so, we can permit each element of X to be an upper bound of the empty chain $\varnothing$. \newline\newline
In all that follows, $X\neq\varnothing$.
\newline\newline
Let $\bar{s}$ be the function form X to P(X), given by $\bar{s}$(x)=$\{ y\in X:y\leq x\}$, and S=ran($\bar{s}$) be partially ordered by inclusion.\newline
Then $\bar{s}$ is a one-to-one function, as, if not, $\exists x,y\in X$ with $x\neq y$, and $\bar{s}(x)=\bar{s}(y)$. But as $\bar{s}(x)=\bar{s}(y)$, we have $x\leq y$ and $y\leq x$. Therefore, x=y, which is a contradiction.\newline
If $\exists x,y\in X$ with $x\leq y$, then $\bar{s}(x)\subset \bar{s}(y)$, as all elements less than x are less than or equal to y. On the other hand, if $\bar{s}(x)\subset \bar{s}(y)$, then $x\leq y$ because $x\in\bar{s}(x)$.\newline
Thus, we have the following Lemma.
\begin{lemma}
A necessary and sufficient condition for $\bar{s}(x)\subset \bar{s}(y)$ is $x\leq y$.
\end{lemma}
This, immediately gives us the next Lemma.
\begin{lemma}
If $\exists x\in X$ such that x is maximal in X, then $\bar{s}(x)$ is maximal in S and vice versa.
\end{lemma}

Let $\mathscr{X}$ be the set of all chains in X, ordered by inclusion. Then, $\mathscr{X}\neq\varnothing$ as singletons and $\varnothing$ in P(X) are chains in X.\newline\newline
Suppose that $\mathscr{Y}\in\mathscr{X}$. By the hypothesis of Zorn's Lemma, $\mathscr{Y}$ has an upper bound in X, say z. Then, $\mathscr{Y}\subset\bar{s}(z)$.\newline\newline
Let $\mathscr{C}$ be a chain in $\mathscr{X}$, and Y be $\cup\{ C:C\in\mathscr{C} \}$.\newline
If $\mathscr{C}=\varnothing$, then $Y=\varnothing$(the empty chain in X), which is in fact in $\mathscr{X}$.
\newline Otherwise, $\exists$x,y$\in Y$, and $\exists C_x,C_y\in \mathscr{C}$ with $x\in C_x$ and $y\in C_y$. Moreover, $\mathscr{C}$ being a chain in $\mathscr{X}$, either $C_x\subset C_y$ or $C_y\subset C_x$.\newline
If $C_x\subset C_y$, then x,y$\in C_y$ and therefore, x and y are comparable.\newline
If $C_y\subset C_y$, then x,y$\in C_x$ and therefore, x and y are comparable.\newline
The conclusion is that Y$\in\mathscr{X}$ for every chain  $\mathscr{C}$ in  $\mathscr{X}$ .\newline
Also, for all $C\in\mathscr{C},C\subset Y$. Hence, Y is an upper bound for $\mathscr{C}$ in $\mathscr{X}$.\newline
Therefore,
\begin{lemma}
Every chain in $\mathscr{X}$ has an upper bound in $\mathscr{X}$, the union of the elements of the chain being one.
\end{lemma}

We consider the case, when $\mathscr{X}$  has a maximal element, say N. As N is a chain in X, it has an upper bound in X, say n. Therefore, $N\cup \{ n\}$ is a chain in X with $N\subset N\cup \{ n\}$. But, as N is maximal in $\mathscr{X}$, $N\cup \{ n\}=N$. So, $n\in N$ and is therefore the greatest element of N.\newline
Hence,
\begin{lemma}
A maximal chain in X has a unique upper bound in $\mathscr{X}$, which turns out to be its greatest element.
\end{lemma} 
Note that $\varnothing$ is not a maximal chain in X, as it is included in every singleton in P(X).\newline\newline 
We continue to study the case, where N is a maximal element of  $\mathscr{X}$.\newline We claim that $\bar{s}(n)$ is a maximal element of S, where n is the upper bound of N.\newline
 Let, if possible, $\bar{s}(n)$ not be a maximal element of S. Then, $\exists t\in X$, such that $\bar{s}(n)\subsetneq\bar{s}(t)$, i.e., $n<t$ by Lemma 1. As $N\subset\bar{s}(n)$, we have $N\subset\bar{s}(t)$ and therefore, t is an upper bound of N.\newline But, by Lemma 4, t=n, which leads to a contradiction.\newline
Therefore, using Lemma 2, we have the following conclusion.
\begin{lemma}
If N is a maximal chain in X with the upper bound n, then $\bar{s}(n)$ is a maximal element of S and n is a maximal element of X.
\end{lemma}

\mybox{gray}{To complete the proof of Zorn's Lemma, it is enough to show that $\mathscr{X}$ has a maximal element.}\newline\newline
It will be more convenient and revealing to consider a general setup of a set $Z\subset P(X)$, satisfying:\newline
1) $\varnothing\in Z$,\newline
2) if $A\in Z$ and $B\subset A$, then $B\in Z$,\newline
3) if $\mathscr{C}$ is a chain in Z, $\cup\{ C:C\in\mathscr{C} \}\in Z$\newline
and partially ordered by inclusion.\newline\newline
\mybox{gray}{Note that $\mathscr{X}$ satisfies these properties, and thus our problem now reduces to proving that there exists a maximal element in Z.}\newline\newline
\newline For $A\in Z$, let $\hat{A}=\{ x\in X: A\cup\{ x\}\in Z\}$. Clearly, $A\subset\hat{A}$.
\newline\newline
Note that $\hat{\varnothing}=\{ x\in X:\{ x\}\in Z\}$.\newline
If $Z=\{\varnothing\}$, there exist no singletons in Z, and hence, $\hat{\varnothing}=\varnothing$. Also, if $\hat{\varnothing}=\varnothing$, then there exist no singletons in Z, and due to property 2) of Z, $Z=\{\varnothing\}$. Therefore, $\hat{\varnothing}=\varnothing$ if and only if $Z=\{\varnothing\}$.\newline
On the other hand, if $Z\neq\{\varnothing\}$. Then $\exists C\in Z$ with $C\neq\varnothing$. As a result, $\exists z\in C$ and so $\{ z\}\subset C$. From Property 2) of Z, $\{ z\}\in Z$. So, $\hat{\varnothing}=\cup\{ C : C\in Z\}\neq\varnothing$. 
\newline\newline
As explained in the description of the Axiom of Choice, there exists a function $f:P(X)\setminus\{\varnothing\}\to X$, such that $f(A)\in A$ for all $A\in dom(f)$.\newline\newline
Define a function $g:Z\to Z$ as follows:\newline
a) If $\hat{A}\setminus A\neq\varnothing$, then g(A)=$A\cup\{ f(\hat{A}\setminus A)\}$. Here, $g(A)\in Z$, as $f(\hat{A}\setminus A)\in\hat{A}\setminus A$.\newline
b) If $\hat{A}\setminus A=\varnothing$, then g(A)=A.\newline
Also, suppose that g(A)=A but $\hat{A}\setminus A\neq\varnothing$. Then $\{ f(\hat{A}\setminus A)\}\in A$ but by definition, $\{ f(\hat{A}\setminus A)\}\in\hat{A}\setminus A$, which leads to a contradiction.\newline
Therefore, $\hat{A}\setminus A=\varnothing$, if and only if g(A)=A.\newline\newline
Clearly $A\subset g(A)$. \newline\newline
Also, note that $g(\varnothing)=f(\hat{\varnothing})$ if and only if $Z\neq\{\varnothing\}$, and $g(\varnothing)=\varnothing$ if and only if $Z=\{\varnothing\}$.\newline\newline
Suppose that $\exists A\in Z$, such that $\hat{A}\setminus A=\varnothing$ and A is not maximal in Z. Then, there exists a $C\in Z$, such that $A\subsetneq C$. Therefore, $\exists x\in C$, such that $x\notin A$. As $C\in Z$, by property 2) of Z, $A\cup\{ x\}\in Z$ and so $x\in \hat{A}$. But, as $x\in\hat{A}$ and $x\notin A$, so $\hat{A}\setminus A\neq\varnothing$, which leads to a contradiction.\newline
Also, if A is maximal in Z, there are no elements in X that can be adjoined to A to create a bigger set present in Z, and so $\hat{A}=A$.\newline
Hence, we have the following Lemma.
\begin{lemma}
A is maximal in Z, if and only if $\hat{A}\setminus A=\varnothing$, if and only if g(A)=A.\newline
\end{lemma}
\mybox{gray}{Our problem now reduces to proving that there exists an A in Z such that g(A)=A.}\newline\newline
Note that, although $\varnothing$ satisfies this requirement if $Z=\{\varnothing\}$, it does not suit our purpose of finding a maximal element in X, as, if $Z=\{\varnothing\}$, then $\mathscr{X}=\{\varnothing\}$, which means that $X=\varnothing$, which is not true due to the discussion in the beginning of the proof. Therefore, in further discussion,  the search for a maximal element in Z will not include $\varnothing$.\newline\newline
 Define a subset J of Z, to be a tower if:\newline
1) $\varnothing\in J$,\newline
2) if $A\in J$, then $g(A)\in J$, and\newline
3) if $\mathscr{C}$ is a chain in J, $\cup\{ C:C\in\mathscr{C} \}\in J$.\newline
Such a J does exist, as Z itself satisfies these properties.\newline\newline
Note that $J=\{\varnothing\}$ if and only if $Z=\{\varnothing\}$. The reason is that, if $Z=\{\varnothing\}$, then the only subset of Z which satisfies the properties of a tower is $\{\varnothing\}$, and if $J=\{\varnothing\}$, then by property 2) of J, $g(\varnothing)\in J$ and so  $g(\varnothing)=\varnothing$ which leads to $Z=\{\varnothing\}$, as explained before.
Therefore, from the discussion above the definition of towers, $J=\{\varnothing\}$ is not allowed.
\newline\newline
Also, the intersection of a non empty family $\{ A_i\}$ of towers is a tower, as:\newline
a) $\varnothing\in A$ for all i, therefore, $\varnothing\in\cap _i A_i$,\newline
b) if $x\in \cap _i A_i$, then $x\in A_i$ for all i, and so by condition 2) in the definition of towers, g(x)$\in A_i$ for\newline all i. Therefore, g(x)$\in\cap _i A_i$,\newline
c) if $\mathscr{C}$ is a chain in $\cap _i A_i$, then it is a chain in $A_i$ for all i and so $\cup\{ C:C\in\mathscr{C}\}\in A_i$ for all i, and therefore $\cup\{ C:C\in\mathscr{C}\}\in \cap _i A_i$.\newline
Therefore,
\begin{lemma}
Let $J_o$ be the intersection of all towers, then $J_o$ is the smallest tower.
\end{lemma}

Let $B\in J_o$, be called comparable, if for all $A\in J_o$, either $A\subset B$ or $B\subset A$.
Comparable sets do exist, as $\varnothing\in J_o$ and for all $A\in J_o$, $\varnothing\in A$.\newline\newline
Consider a particular comparable element $\mathfrak{C}$.\newline
Suppose that $A\subsetneq  \mathfrak{C}$. As  $\mathfrak{C}$ is comparable and g(A)$\in J_o$ (as $A\in J_o$), therefore $g(A)\subset  \mathfrak{C}$ or $ \mathfrak{C}\subsetneq g(A)$. But $ \mathfrak{C}\not \subsetneq g(A)$ as, if $ \mathfrak{C}\subsetneq g(A)$, then $A\subsetneq  \mathfrak{C}\subsetneq g(A)$, but g(A) has at most one more element than A.\newline
Therefore,\newline
if $A\subsetneq  \mathfrak{C}$, then $g(A)\subset  \mathfrak{C}$.
\newline\newline
For the particular comparable element $\mathfrak{C}$,\newline
let U=$\{ B\in J_o:B\subset \mathfrak{C}$ or $g(\mathfrak{C})\subset B\}$, partially ordered by inclusion.\newline\newline
Then, $\varnothing\in U$ as $\varnothing\in J_o$ and $\varnothing\subset \mathfrak{C}$.\newline\newline
Also, if $A\in U$, then $g(A)\in U$ due to the following:\newline
As $A\in J_o$, $g(A)\in J_o$.\newline
Case 1: If $A\subsetneq \mathfrak{C}$, then $g(A)\subset \mathfrak{C}$, as proved above. Therefore, $g(A)\in U$.\newline
Case 2: If $A=\mathfrak{C}$, then $g(A)=g(\mathfrak{C})$ and so $g(\mathfrak{C})\subset g(A)$. Therefore, $g(A)\in U$.\newline
Case 3: If $g(\mathfrak{C})\subset A$, then as $A\subset g(A)$, $g(\mathfrak{C})\subset g(A)$. Therefore, $g(A)\in U$.\newline\newline
Consider a chain $\mathscr{C}$ in U.\newline
For any $E\in\mathscr{C}$, $E\in U$ and so E can be of two types:\newline
Type 1: $E\subset\mathfrak{C}$.\newline
Type 2: $g(\mathfrak{C})\subset E$.\newline
We note that the only element of $\mathscr{C}$ of both Type 1 and Type 2, if any, is $\mathfrak{C}$, and that too, if and only if $\mathfrak{C}$ is maximal in Z, as $g(\mathfrak{C})\subset E\subset\mathfrak{C}$ implies that $E=\mathfrak{C}=g(\mathfrak{C})$.\newline
Let Y=$\cup\{ C:C\in\mathscr{C}\}$. There are two possibilities for Y:\newline
Possibility 1: All elements of $\mathscr{C}$ are of Type 1. Then, all elements of Y are in $\mathfrak{C}$ and so $Y\subset \mathfrak{C}$. \newline
Possibility 2: Atleast one of the Cs in $\mathscr{C}$, say E, is of Type 2. Then, $g(\mathfrak{C})\subset E$. As $E\subset Y$, therefore $g(\mathfrak{C})\subset Y$.\newline
As a result, $Y\in U$.\newline\newline
From the previous three paragraphs, U is a tower included in $J_o$.\newline
But due to Lemma 7, U=$J_o$. Therefore, for all $A\in J_o$, $A\in U$ and so $A\subset \mathfrak{C}$ or $g(\mathfrak{C})\subset A$.\newline
As $A\subset\mathfrak{C}$ implies that $A\subset g(\mathfrak{C})$, therefore for all $A\in J_o$, $A\subset g(\mathfrak{C})$ or $g(\mathfrak{C})\subset A$.\newline
As a result,
\begin{lemma}
If $\mathfrak{C}$ is comparable, then $g(\mathfrak{C})$ is also comparable.
\end{lemma}

Now let $C_o$ be the set of all comparable sets. Clearly, $C_o$ is a chain in $J_o$, as comparable sets are comparable with each other.\newline\newline
1. $\varnothing\in C_o$.\newline
2. From Lemma 8, if $C\in C_o$, then $g(C)\in C_o$.\newline
3. Consider a chain $\mathscr{C}$ in $C_o$.\newline
If $E\in\mathscr{C}$, then for a particular $A\in J_o$, E can be of two types:\newline
Type 1: $A\subset E$.\newline
Type 2: $E\subset A$.\newline
Note that the only element of $\mathscr{C}$ of both Type 1 and Type 2, if any, is A itself.\newline
Let Y=$\cup\{ C:C\in\mathscr{C}\}$. There are two possibilities for Y:\newline
Possibility 1: All elements of $\mathscr{C}$ are of Type 2. Then, all elements of Y are in A, i.e., $Y\subset A$.\newline
Possibility 2: Atleast one of the Cs in $\mathscr{C}$, say E, is of Type 1. Then, $A\subset E$. As, $E\subset Y$, therefore $A\subset Y$.\newline
Therefore, Y is comparable with A. As A was an arbitrary element of $J_o$, Y is a comparable set, and hence, $Y\in C_o$.
\newline\newline

From 1), 2) and 3) in the previous paragraph, we infer that $C_o$ is a tower. But by Lemma 7, $C_o=J_o$.
Hence, the important result:
\begin{lemma}
$J_o$ is a chain in $J_o$.\newline
\end{lemma}
Now, let A=$\cup\{ C:C\in J_o\}$.\newline
Due to Lemma 9 and property 3) of towers, $A\in J_o$.\newline
By property 2) of towers, as $A\in J_o$, therefore $g(A)\in J_o$.\newline
Also, as A is the union of all sets in $J_o$, $g(A)\subset A$.\newline
But as $A\subset g(A)$ always, therefore \textbf{A=g(A)}.\newline\newline
Thus, we have obtained the desired result.
\end{proof} 

\end{document}